\documentclass[final,leqno]{siamltex}

\usepackage{amsmath, amssymb}
\usepackage[table]{xcolor}
\usepackage{placeins}
\usepackage{url}

\usepackage{tikz,pgfplots}
\usepackage{graphicx}
\usepackage{subfigure}
\usepackage{placeins}
\usepackage[left=0.5in, right=0.5in, top=1.5in, bottom=1.5in]{geometry}
\newtheorem{claim}[theorem]{Claim}
\newtheorem{defn}[theorem]{Defn}

\usepackage[numbers,sort,compress]{natbib}

\title{Generalized Rybicki Press algorithm}
\subtitle{$\mathcal{O}(N)$ direct solver and determinant computation of exponential covariance and general semi-separable matrices}

\author{Sivaram Ambikasaran}

\definecolor{mycyan}{rgb}{0,0.95,0.95}

\pgfplotsset{compat=newest}

\begin{document}

\maketitle

\begin{abstract}
This article discusses a more general and numerically stable Rybicki Press algorithm, which enables inverting and computing determinants of covariance matrices, whose elements are sums of exponentials. The algorithm is true in exact arithmetic and relies on introducing new variables and corresponding equations, thereby converting the matrix into a banded matrix of larger size. Linear complexity banded algorithms for solving linear systems and computing determinants on the larger matrix enable linear complexity algorithms for the initial semi-separable matrix as well. Benchmarks provided illustrate the linear scaling of the algorithm.
\end{abstract}



\pagestyle{myheadings}
\thispagestyle{plain}
\markboth{$\mathcal{O}(N)$ direct solver and determinant computation of general semi-separable matrices}{Generalized Rybicki Press algorithm}

\section{Introduction}
Large dense covariance matrices arise in a wide range of applications in computational statistics and data analysis. Storing and performing numerical computations on such large dense matrices is computationally intractable. However, most of these large dense matrices are structured (either in exact arithmetic or finite arithmetic), which can be exploited to construct fast algorithms. One such class of data sparse matrices are semi-separable matrices, which have raised a lot of interest and have been studied in detail across a wide range of applications including integral equations~\cite{asplund1959finite,gohberg1984non,gesztesy2003modified} and computational statistics~\cite{roy1956inverting,roy1960evaluation,mustafi1967inverse,uppuluri1969inverse, greenberg1959matrix}. For a detailed bibliography on semi-separable matrices, the reader is referred to Vandebril et al.~\cite{vandebril2005bibliography}. Throughout the literature, there are slightly different definitions of semi-separable matrices. In this article, we will be working with the following definition:
\begin{defn}
$A \in \mathbb{R}^{N \times N}$ is termed a semi-separable matrix with semi-separable rank $p$, if it can be written as
\begin{align}
A & = D + \text{triu}(B_p) + \text{tril}(C_p)
\end{align}
where $D$ is a diagonal matrix, $B_p, C_p$ are rank $p$ matrices, $\text{triu}(B_p)$ denotes the upper triangular part of $B_p$ and $\text{tril}(C_p)$ denotes the lower triangular part of $C_p$.
\end{defn}

Fast algorithms for solving semi-separable linear systems exists and the reader is referred to some of these references~\cite{eidelman2002modification,van2004two,eidelman1997inversion,gohberg1985linear, jain2010n,chandrasekaran2003fast} and the references therein. In this article, we propose a new $\mathcal{O}(N)$ direct solver and determinant computation for semi-separable matrices.

The main contributions of this article include:

\begin{itemize}
\item
A new $\mathcal{O}(N)$ direct solver for semi-separable matrices is obtained by embedding the semi-separable matrix into a larger banded matrix.
\item
The determinant of these semi-separable matrix is shown to equal to the determinant of the larger banded matrix, thereby enabling computing determinants of these semi-separable matrices at a computational cost of $\mathcal{O}(N)$. This is the first algorithm for computing the determinants for a general semi-separable matrix.
\item
A numerically stable generalized Rybicki Press algorithm is derived using these ideas. To be specific, fast, stable, direct algorithms are derived for solving and computing determinants (both scaling as $\mathcal{O}(N)$) for covariance matrices of the form:
\begin{align}
A_{ij} = \sum_{l=1}^p \alpha_l \exp\left(-\beta_l \vert t_i - t_j \vert\right)
\label{eqn_covariance_matrix}
\end{align}
where $i,j \in \{1,2,\ldots,n\}$, the points $t_i$ are distinct and are distributed on an interval. The covariance matrix in Equation~\eqref{eqn_covariance_matrix} is frequently encountered in computational statistics in the context of Continuous time AutoRegressive-Moving-Average (abbreviated as CARMA) models~\cite{brockwell2002introduction, brockwell2001levy, brockwell1994continuous}.
\item
Another advantage of this algorithm from a practical view-point is that the algorithm relies only on sparse linear algebra and thereby can easily use the existing mature sparse linear algebra libraries.
\end{itemize}

The algorithm discussed in this article has been implemented in C++ and the implementation is made available at \url{https://github.com/sivaramambikasaran/ESS}~\cite{ambikasaran2014ESS} under the license provided by New York University.

\textbf{Acknowledgements}: The author would like to thank Christopher S. Kochanek for initiating the conversation on generalized Rybicki Press algorithm and David W. Hogg for putting in touch with Christopher S. Kochanek. The author would also like to thank the anonymous referee for his careful, detailed review and insightful comments. The research was supported in part by the NYU-AIG Partnership on Innovation for Global Resilience under grant number A2014-005. The author was also supported in part by the Applied Mathematical Sciences Program of the U.S. Department of Energy under Contract DEFGO288ER25053, Office of the Assistant Secretary of Defense for Research and Engineering and AFOSR under NSSEFF Program Award FA9550-10-1-0180.

\section{Sparse embedding of semi-separable matrix with semi-separable rank $1$}
To motivate the general idea, we will first look at the sparse embedding for a $4 \times 4$ semi-separable matrix, whose semi-separable rank is $1$. The matrix $A$ is as shown in Equation~\eqref{eqn_first_example}.

\begin{align}
A & =
\begin{bmatrix}
a_{11} & u_1v_2 & u_1v_3 & u_1v_4\\
u_1v_2 & a_{22} & u_2v_3 & u_2v_4\\
u_1v_3 & u_2v_3 & a_{33} & u_3v_4\\
u_1v_4 & u_2v_4 & u_3v_4 & a_{44}\\
\end{bmatrix}
\label{eqn_first_example}
\end{align}
And the corresponding linear system is $Ax=b$, where $b=\begin{bmatrix} b_1& b_2& b_3& b_4 \end{bmatrix}^T$

Introduce the following variables:
\begin{align}
r_4 & = v_4x_4\\
r_3 & = v_3x_3 + r_4\\
r_2 & = v_2x_2 + r_3
\end{align}
\begin{align}
l_1 & = u_1x_1\\
l_2 & = u_2x_2 + l_1\\
l_3 & = u_3x_3 + l_2
\end{align}

Introducing the variables the linear system $Ax=b$ is now of the form
\begin{align}
a_{11}x_1 + u_1r_2 & = b_1\\
v_2 l_1 + a_{22}x_2 + u_2r_3 & = b_2\\
v_3 l_2 + a_{33}x_3 + u_3r_4 & = b_3\\
v_4 l_3 + a_{44}x_4 & = b_4
\end{align}

\begin{figure}[!htbp]
\begin{center}
\includegraphics[scale=0.5]{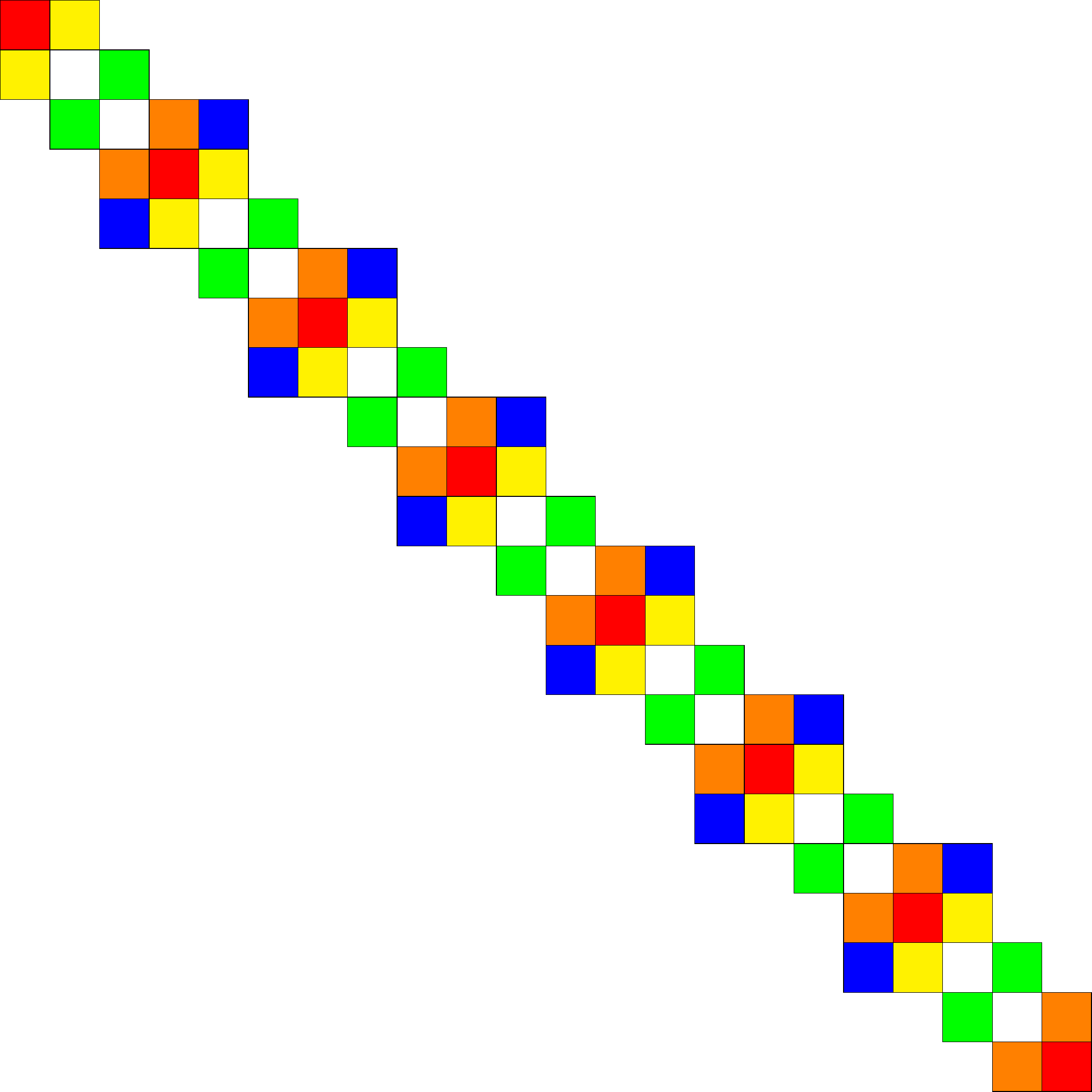}
\end{center}
\caption{Pictorial description of the extended sparse matrix obtained from a rank $1$ semi-separable matrix where $N=8$. The color code is as shown below.}
\begin{center}
\includegraphics[scale=1]{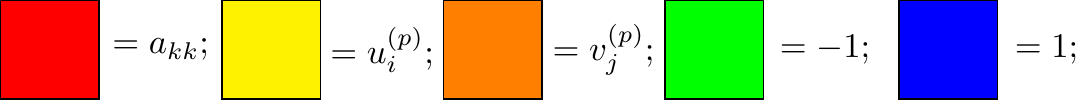}
\end{center}
\label{fig_rank1_semiseparable}
\end{figure}

The extended linear system (after appropriate ordering of equations and unknowns) is then of the form
\begin{align}
\begin{bmatrix}
a_{11} & u_1 & 0 & 0 & 0 & 0 & 0 & 0 & 0 & 0\\
u_1 & 0 & -1 & 0 & 0 & 0 & 0 & 0 & 0 & 0\\
0 & -1 & 0 & v_2 & 1 & 0 & 0 & 0 & 0 & 0\\
0 & 0 & v_2 & a_{22} & u_2 & 0 & 0 & 0 & 0 & 0\\
0 & 0 & 1 & u_2 & 0 & -1 & 0 & 0 & 0 & 0\\
0 & 0 & 0 & 0 & -1 & 0 & v_3 & 1 & 0 & 0\\
0 & 0 & 0 & 0 & 0 & v_3 & a_{33} & u_3 & 0 & 0\\
0 & 0 & 0 & 0 & 0 & 1 & u_{3} & 0 & -1 & 0\\
0 & 0 & 0 & 0 & 0 & 0 & 0 & -1 & 0 & v_4\\
0 & 0 & 0 & 0 & 0 & 0 & 0 & 0 & v_4 & a_{44}\\
\end{bmatrix}
\begin{bmatrix}
x_1\\
r_2\\
l_1\\
x_2\\
r_3\\
l_2\\
x_3\\
r_4\\
l_3\\
x_4
\end{bmatrix}
& =
\begin{bmatrix}
b_1\\
0\\
0\\
b_2\\
0\\
0\\
b_3\\
0\\
0\\
b_4
\end{bmatrix}
\label{eqn_extended_example}
\end{align}
Note that Equation~\eqref{eqn_extended_example} is a banded matrix of bandwidth $2$ and has a sparsity structure even within the band. In general, let $A$ be an $N \times N$ semi-separable matrix, with the semi-separability rank $1$ as written in Equation~\eqref{eqn_generic_semiseparable_rank_1}.
\begin{align}
A(i,j) & =
\begin{cases}
a_{ii} & \text{ if }i=j\\
u_j v_i & \text{ if }i> j\\
u_i v_j & \text{ if }i< j\\
\end{cases}
\label{eqn_generic_semiseparable_rank_1}
\end{align}
where $i,j \in \{1,2,\ldots,N\}$. One would then need to add the variables $r_2,r_3,\ldots,r_N$ and $l_1,l_2,\ldots,l_{N-1}$, where $r_N = v_N x_N$, $l_1 = u_1x_1$ and
\begin{align}
r_k & = v_kx_k + r_{k+1}\\
l_k & = u_kx_k + l_{k-1}
\end{align}
where $k \in \{2,\ldots,N-1\}$. Hence, we have a total of $3N-2$ variables and $3N-2$ equations. Therefore, the extended matrix will be a $(3N-2) \times (3N-2)$ banded matrix, whose bandwidth is $2$. This is illustrated pictorially for a $8 \times 8$ matrix in Figure~\ref{fig_rank1_semiseparable}.

\section{Sparse embedding of a general semi-separable matrix}
\begin{figure}[!htbp]
\begin{center}
\includegraphics[scale=0.225]{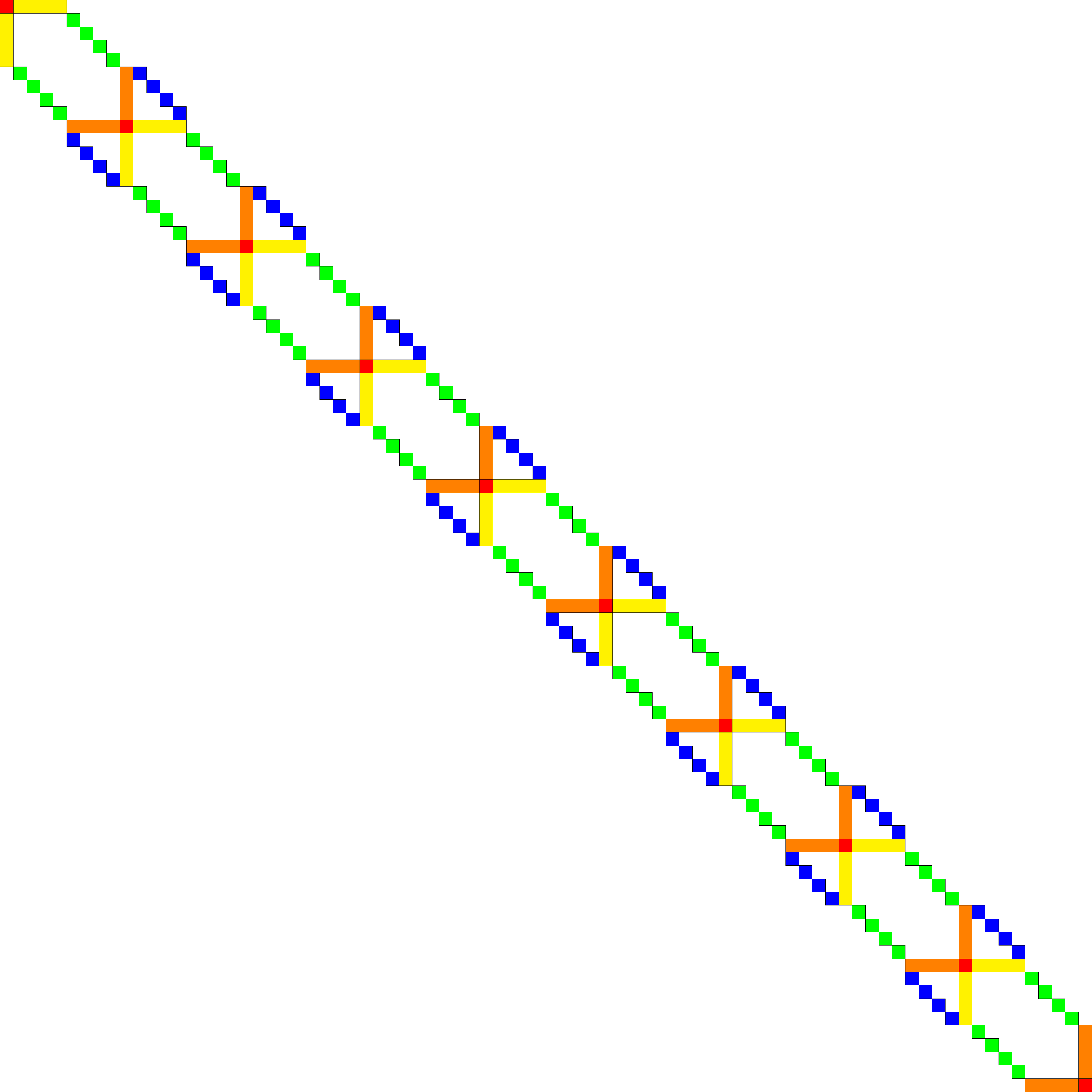}
\end{center}
\caption{Pictorial description of the extended sparse matrix where $N=10$ and $p=4$. The color code is as shown below.}
\begin{center}
\includegraphics[scale=1]{./images/colorcode.pdf}
\end{center}
\label{fig_rankm_semiseparable}
\end{figure}

Let $A$ be a $N \times N$ matrix, whose semi-separable rank is $p$, i.e., we have
\begin{align}
A(i,j) & =
\begin{cases}
a_{ii} & \text{ if }i=j\\
\displaystyle\sum_{l=1}^p u_j^{(l)} v_i^{(l)} & \text{ if }i>j\\
\displaystyle\sum_{l=1}^p u_i^{(l)} v_j^{(l)} & \text{ if }i<j\\
\end{cases}
\end{align}
where $i,j \in \{1,2,\ldots,N\}$. We then add the following variables $r_2^{(p)},r_3^{(p)},\ldots,r_N^{(p)}$ and $l_1^{(p)},l_2^{(p)},\ldots,l_{N-1}^{(p)}$ as before. However, not surprisingly, these new variables $r_i^{(p)}$'s and $l_j^{(p)}$'s will be vectors of length $p$. Let $U_k^{(p)} = \begin{bmatrix} u_k^{(1)} & u_k^{(2)} & u_k^{(3)} & \cdots & u_k^{(p)}\end{bmatrix}$ and $V_k^{(p)} = \begin{bmatrix} v_k^{(1)} & v_k^{(2)} & v_k^{(3)} & \cdots & v_k^{(p)}\end{bmatrix}$. We then have the following relations for the additional vector variables.
\begin{align}
r_N^{(p)} = V_N^Tx_N\\
l_1^{(p)} = U_1^Tx_1
\end{align} and
\begin{align}
r_k^{(p)} & = V_k^T x_k + r_{k+1}^{(p)}\\
l_k^{(p)} & = U_k^T x_k + l_{k-1}^{(p)}
\end{align}
where $k \in \{2,\ldots,N-1\}$. Hence, we now have $(2p+1)N-2p$ variables (this includes the $N$ $x_k$'s, $N-1$ vector variables $r_k^{(p)}$ and $l_k^{(p)}$ of length $p$) and $(2p+1)N-2p$ equations relating them. Therefore, we end up with a $((2p+1)N-2p) \times ((2p+1)N-2p)$ extended sparse matrix, whose bandwidth is $(2p+1)$. This is illustrated in Figure~\ref{fig_rankm_semiseparable} for $10 \times 10$ semi-separable matrix, whose semi-separable rank is $4$.

The computational complexity of the algorithm clearly scales as $\mathcal{O}(N)$, since the extended sparse matrix has a bandwidth of $\mathcal{O}(p)$ and the matrix of size $\mathcal{O}(pN) \times \mathcal{O}(pN)$. It is also possible to analyze the scaling with respect to the semi-separable rank $p$, though this is of little practical relevance since $p = \mathcal{O}(1)$ for most interesting semi-separable matrices. A detailed analysis shows that the computational complexity of the algorithm is $\mathcal{O}(p^2N)$. Numerical benchmarks presented in Section~\ref{section_nb} validate the scaling of the algorithm.


\section{Determinant of extended sparse matrix}

\begin{claim}
The determinant of the extended sparse matrix is the same as the determinant of the original dense matrix up to a sign.
\end{claim}

The extended system, denoted by $A_{ex}$ on appropriate reordering of rows and columns can be written as
\begin{align}
P_1A_{ex} P_2
\begin{bmatrix}
l_1^{(p)}\\
l_2^{(p)}\\
\vdots\\
l_{N-1}^{(p)}\\
r_2^{(p)}\\
r_3^{(p)}\\
\vdots\\
r_N^{(p)}\\
x_1\\
x_2\\
\vdots\\
x_N
\end{bmatrix} & =
\begin{bmatrix}
L_{\Delta} & 0 & U_a\\
0 & U_{\Delta} & V_a\\
V_b & U_b & D
\end{bmatrix}
\begin{bmatrix}
l_1^{(p)}\\
l_2^{(p)}\\
\vdots\\
l_{N-1}^{(p)}\\
r_2^{(p)}\\
r_3^{(p)}\\
\vdots\\
r_N^{(p)}\\
x_1\\
x_2\\
\vdots\\
x_N
\end{bmatrix}
\end{align}
where $P_1$, $P_2$ are permutation matrices, the matrix $L_{\Delta}$ is a highly sparse lower-triangular matrix with $1$'s on the diagonal and $-1$'s at a few places in the lower-triangular part (the precise location is unimportant for determinant computations as we will see later), the matrix $U_{\Delta}$ is a highly sparse upper-triangular matrix with $1$'s on the diagonal and $-1$'s at a few places in the upper-triangular part and $D$ is a diagonal matrix with $D_{ii} = a_{ii}$. The first set of rows, i.e., $\begin{bmatrix} L_{\Delta} & 0 & U_a\end{bmatrix}$, correspond to adding the variables $l_k^{(p)}$, i.e., $l_k^{(p)} = U_k^Tx_k + l_{k-1}^{(p)}$. The next set of rows, i.e., $\begin{bmatrix} 0 & U_{\Delta} & V_a\end{bmatrix}$, correspond to adding the variables $r_k^{(p)}$, i.e., $r_k^{(p)} = V_k^Tx_k + r_{k+1}^{(p)}$. The last set of rows, i.e., $\begin{bmatrix} V_b & U_b & D\end{bmatrix}$, correspond to the initial set of equations with the $l_k^{(p)}$'s and $r_k^{(p)}$'s introduced.
We then have
\begin{align}
\det(P_1A_{ex}P_2) & =
\underbrace{
\det\left(
\begin{bmatrix}
L_{\Delta} & 0 & U_a\\
0 & U_{\Delta} & V_a\\
V_b & U_b & D
\end{bmatrix}
\right) =
\det\left(
\begin{bmatrix}
L_{\Delta} & 0\\
0 & U_{\Delta}
\end{bmatrix}
\right)
\det\left(D-
\begin{bmatrix}
V_b & U_b
\end{bmatrix}
\begin{bmatrix}
L_{\Delta} & 0\\
0 & U_{\Delta}
\end{bmatrix}^{-1}
\begin{bmatrix}
U_a\\
V_a
\end{bmatrix}
\right)
}_{\text{Block determinant formula}}
\end{align}
Now note that $\det(L_{\Delta}) = 1 = \det(U_{\Delta})$, due to the fact that $L_{\Delta}$ and $U_{\Delta}$ are triangular matrices with $1$'s on the diagonal. Hence,
\begin{align}
\det\left(
\begin{bmatrix}
L_{\Delta} & 0\\
0 & U_{\Delta}
\end{bmatrix}
\right) & = \det(L_{\Delta}) \det(U_{\Delta}) = 1 \times 1 = 1
\end{align}
Further, note that the matrix $D-
\begin{bmatrix}
V_b & U_b
\end{bmatrix}
\begin{bmatrix}
L_{\Delta} & 0\\
0 & U_{\Delta}
\end{bmatrix}^{-1}
\begin{bmatrix}
U_a\\
V_a
\end{bmatrix}
$ is the Schur complement obtained by eliminating the variables $l_i^{(p)}$, $r_i^{(p)}$ and hence is the initial dense matrix $A$ we began with, i.e.,
\begin{align}
D-
\begin{bmatrix}
V_b & U_b
\end{bmatrix}
\begin{bmatrix}
L_{\Delta} & 0\\
0 & U_{\Delta}
\end{bmatrix}^{-1}
\begin{bmatrix}
U_a\\
V_a
\end{bmatrix} = A
\end{align}
Hence, we have
\begin{align}
\det(P_1A_{ex}P_2) & = \det\left(D-
\begin{bmatrix}
V_b & U_b
\end{bmatrix}
\begin{bmatrix}
L_{\Delta} & 0\\
0 & U_{\Delta}
\end{bmatrix}^{-1}
\begin{bmatrix}
U_a\\
V_a
\end{bmatrix}\right) = \det(A)
\end{align}
which gives us that
\begin{align}
\det(A_{ex}) = \pm \det(A)
\end{align}
where the ambiguity in the sign arises due to the determinant of the permutation matrices.
\section{Reinterpretation of Rybicki Press algorithm in terms of sparse embedding}

We will first naively reinterpret the Rybicki Press algorithm in terms of the extended sparse matrix algebra. Recall that the Rybicki Press algorithm~\cite{rybicki1995class} inverts a correlation matrix $A$ given by Equation~\eqref{eqn_Rybicki_Press_matrix}.
\begin{align}
A(i,j) & = \exp\left(-\beta \lvert t_i - t_j \rvert\right)
\label{eqn_Rybicki_Press_matrix}
\end{align}
where $t_i$'s lies on an interval and are monotone. The original Rybicki Press algorithm relies on the fact that the inverse of $A$ happens to be a tridiagonal matrix. The key ingredient of their algorithm is the following property of exponentials:
\begin{align}
\exp\left(\beta (t_i-t_j) \right)\exp\left(\beta (t_j-t_k) \right) = \exp\left(\beta (t_i-t_k)\right)
\end{align}
In our sparse interpretation as well, we will use this property to recognize that the matrix $A$ is a semi-separable matrix, whose semi-separable rank is $1$. This can be seen by setting $u_k = \exp(\beta t_k)$ and $v_k = \exp(-\beta t_k)$. This then gives us ($i<j$) that $A(i,j) = u_i v_j = \exp(\beta t_i) \exp(-\beta t_j) = \exp(\beta(t_i-t_j))$ and similarly for $i>j$. This shows that the matrix $A$ is semi-separable with semi-separable rank $1$. Hence, we can mimic the same approach as in the earlier sections to obtain an $\mathcal{O}(N)$ algorithm. However, there is an issue that needs to be addressed from a numerical perspective. If the $t_i$'s are spread over a large interval, then $u_i$ is exponentially large, while $v_i$ is exponentially small, and hence embedding into a sparse matrix as such could prove to be a catastrophic leading to underflow and overflow of the relevant entries. This issue though can be circumvented by a suitable analytic preconditioning, by an appropriate change of variables. This is illustrated for a $4 \times 4$ linear system. We will use the notation $t_{ij}$ to denote $\lvert t_i - t_j\lvert$. The linear equation is
\begin{align}
\begin{bmatrix}
1 & \exp(-\beta t_{12}) & \exp(-\beta t_{13}) & \exp(-\beta t_{14})\\
\exp(-\beta t_{12}) & 1 & \exp(-\beta t_{23}) & \exp(-\beta t_{24})\\
\exp(-\beta t_{13}) & \exp(-\beta t_{23}) & 1 & \exp(-\beta t_{34})\\
\exp(-\beta t_{14}) & \exp(-\beta t_{24}) & \exp(-\beta t_{34}) & 1\\
\end{bmatrix}
\begin{bmatrix}
x_1\\x_2\\x_3\\x_4
\end{bmatrix}
=
\begin{bmatrix}
b_1\\b_2\\b_3\\b_4
\end{bmatrix}
\end{align}
Now lets introduce the additional variables as follows:
\begin{align}
r_4 & = x_4\\
r_3 & = x_3 + \exp(-\beta t_{34})r_4\\
r_2 & = x_2 + \exp(-\beta t_{23})r_3
\end{align}
\begin{align}
l_2 & = x_1\exp(-\beta t_{12})\\
l_3 & = (x_2+l_2)\exp(-\beta t_{23})\\
l_4 & = (x_3+l_3)\exp(-\beta t_{34})
\end{align}
The equations then become
\begin{align}
x_1 + \exp(-\beta t_{12}) r_2 & = b_1\\
l_2 + x_2 + \exp(-\beta t_{23}) r_3 & = b_2\\
l_3 + x_3 + \exp(-\beta t_{24}) r_4 & = b_3\\
l_4 + x_4 & = b_4
\end{align}
Embedding this in an extended sparse matrix, we obtain
\begin{align}
\begin{bmatrix}
1 & \exp(-\beta t_{12}) & 0 & 0 & 0 & 0 & 0 & 0 & 0 & 0\\
\exp(-\beta t_{12}) & 0 & -1 & 0 & 0 & 0 & 0 & 0 & 0 & 0\\
0 & -1 & 0 & 1 & \exp(-\beta t_{23}) & 0 & 0 & 0 & 0 & 0\\
0 & 0 & 1 & 1 & \exp(-\beta t_{23}) & 0 & 0 & 0 & 0 & 0\\
0 & 0 & \exp(-\beta t_{23}) & \exp(-\beta t_{23}) & 0 & -1 & 0 & 0 & 0 & 0\\
0 & 0 & 0 & 0 & -1 & 0 & 1 & \exp(-\beta t_{34}) & 0 & 0\\
0 & 0 & 0 & 0 & 0 & 1 & 1 & \exp(-\beta t_{34}) & 0 & 0\\
0 & 0 & 0 & 0 & 0 & \exp(-\beta t_{34}) & \exp(-\beta t_{34}) & 0 & -1 & 0\\
0 & 0 & 0 & 0 & 0 & 0 & 0 & -1 & 0 & 1\\
0 & 0 & 0 & 0 & 0 & 0 & 0 & 0 & 1 & 1\\
\end{bmatrix}
\begin{bmatrix}
x_1\\ r_2\\ l_1\\ x_2\\ r_3\\ l_2\\ x_3\\ r_4\\ l_3\\ x_4
\end{bmatrix}
=
\begin{bmatrix}
b_1\\ 0\\ 0\\ b_2\\ 0\\ 0\\ b_3\\ 0\\ 0\\ b_4
\end{bmatrix}
\end{align}

Note that the sparsity pattern of the matrix is the same as before, which is to be expected, since all we have done essentially is to scale elements appropriately and hence the zero fill-ins remain the same.

\section{Numerically stable generalized Rybicki Press}
The same idea carries over the generalized Rybicki Press algorithm, i.e., if we consider a CARMA($p$,$q$) process which has the covariance matrix given by
\begin{align}
K(r) & =
\begin{cases}
d & \text{if } r=0\\
\displaystyle \sum_{l=1}^p \alpha_l(p,q) \exp(-\beta_l r) & \text{if }r >0
\end{cases}
\end{align}
then it immediately follows that that the matrix is semi-separable with semi-separable rank being $p$. To avoid numerical overflow and underflow, as shown in the previous section, appropriate sets of variables need to be introduced. Let
\begin{align}
\alpha & = \begin{bmatrix} \alpha_1 & \alpha_2 & \cdots & \alpha_p\end{bmatrix}^T
\end{align}
and
\begin{align}
\gamma_k & = \begin{bmatrix}\exp(-\beta_1 t_{k,k+1})& \exp(-\beta_2 t_{k,k+1}) & \cdots& \exp(-\beta_p t_{k,k+1}) \end{bmatrix}^T
\end{align}
Now introduce the variables
\begin{align}
r_k & = \alpha x_k + D^{(k,k+1)} r_{k+1}
\label{eqn_defn_rk}
\end{align}
\begin{align}
l_{k} & = \gamma_{k-1} x_{k-1} + D^{(k-1,k)}l_{k-1}
\label{eqn_defn_lk}
\end{align}
where $k \in \{2,3,\ldots,N\}$, with $r_{N+1} = l_1 = 0$ and $D^{(k,k+1)}$ is a $p \times p$ diagonal matrix, with its diagonal being $\gamma_k$. The initial equations become
\begin{align}
\alpha^T l_k + dx_k + \gamma_k^T r_{k+1} = b_k
\label{eqn_defn_xk}
\end{align}
where $k \in \{1,2,\ldots,N\}$. Now form the extended sparse matrix using the variables $x_k,l_k$ and $r_k$, with the equations being Equations~\eqref{eqn_defn_rk},~\eqref{eqn_defn_lk},~\eqref{eqn_defn_xk}. The sparsity pattern of the extended sparse matrix is the same and hence the computational complexity scales as $\mathcal{O}(N)$.

\section{Numerical benchmarks}
\label{section_nb}
We present a few numerical benchmarks illustrating the scaling of the algorithm and the error. In all these benchmarks, the semi-separable matrix is of the form
\begin{align}
A(i,j) =
\begin{cases}
d & \text{if }i=j\\
\displaystyle \sum_{l=1}^p \alpha_l \exp\left(-\beta_l \lvert t_i-t_j\rvert\right) & \text{if }i \neq j
\end{cases}
\end{align}
where the $t_i$'s lie on a on-dimensional manifold and are sorted in increasing fashion. Apart from the time taken for the assembly, factorization and solve, the infinity norm of the residual, i.e., $\Vert Ax-b \Vert_{\infty}$ and the relative error in the log determinant are also presented. For the purposes of benchmark, $t_i$'s are chosen at random from the interval $[0,20]$ and then sorted; $\alpha_l$'s, $\beta_l$'s are chosen at random from the interval $[0,2]$; and $d$ is set equal to $1+\displaystyle \sum_{l=1}^p \alpha_l$. Throughout the benchmarks the original dense matrix will be referred to as $A$, while the corresponding extended sparse matrix will be referred to as $A_{ex}$.

The extended sparse linear system, i.e., $A_{ex}x_{ex} = b_{ex}$, is solved using the sparse LU factorization (SparseLU) in Eigen~\cite{eigenweb}. This relies on the sequential SuperLU package~\cite{demmel1999supernodal, demmel2011superlu, li2005overview}, which performs sparse LU decomposition with partial pivoting. The preordering of the unknowns is performed using the COLAMD method~\cite{davis2004algorithm}. It is to be noted that despite the preordering and partial pivoting, which inturn affects the banded structure, the computational cost as shown in Figures~\ref{figure_benchmark1},~\ref{figure_benchmark2} for the extended sparse system scales linearly in the number of unknowns. The extended sparse matrix is stored using a triplet list in Eigen~\cite{eigenweb}, which internally converts it into compressed column/row storage format. The exact implementation can be found at \url{https://github.com/sivaramambikasaran/ESS}~\cite{ambikasaran2014ESS}.

\subsection{Benchmark $1$} 
In this benchmark, we illustrate the linear scaling of the algorithm with the number of unknowns $N$ for different choices of $p$. The solution obtained using the sparse LU factorization is compared with the partial pivoted LU algorithm (PartialPivLU) in Eigen~\cite{eigenweb}, which is used to solve the initial dense linear system $Ax=b$. Table~\ref{table_N_scaling} shows the scaling of the algorithm and the maximum error in the residual for a fixed semi-separable rank of $p=5$.

\begin{table}[!htbp]
\rowcolors{1}{gray!30}{white}
\caption{Scaling of the algorithm with system size $N$ for a fixed semi-separable rank $p=5$. The time taken is reported in milliseconds.}
\begin{center}
\begin{tabular}{|c|c|c|c|c|c|c|c|c|c|}
\hline
System size& \multicolumn{6}{c|}{Time taken in milliseconds} & \multicolumn{2}{c|}{Error in residual} & Error in log-det\\
\hline
$N$ & \multicolumn{2}{c|}{Assembly} & \multicolumn{2}{c|}{Factorize} & \multicolumn{2}{c|}{Solve} & \multicolumn{2}{c|}{measured in $\Vert \cdot \Vert_\infty$} & $\dfrac{\log(\vert A_{ex}\vert/\vert A \vert)}{\log(\vert A \vert)}$\\
\hline
 & Usual & Fast & Usual & Fast & Usual & Fast & Usual & Fast & \\
\hline
$500$ & $15.5$ & $1.15$ & $12$ & $8$ & $0.233$ & $1.36$ & $2 \times 10^{-14}$ & $2.2 \times 10^{-15}$ & $1.08 \times 10^{-15}$\\
\hline
$1000$ & $49.2$ & $1.73$ & $91.6$ & $15.5$ & $0.862$ & $2.02$ & $4 \times 10^{-14}$ & $3.8 \times 10^{-15}$ & $1.46 \times 10^{-15}$\\
\hline
$2000$ & $188$ & $3.28$ & $643$ & $30.8$ & $2.80$ & $4.41$ & $9 \times 10^{-14}$ & $5.6 \times 10^{-15}$ & $1.67 \times 10^{-15}$\\
\hline
$5000$ & $1150$ & $9.11$ & $9360$ & $83.1$ & $14.4$ & $10.4$ & $2 \times 10^{-13}$ & $6.4 \times 10^{-15}$ & $5.44 \times 10^{-16}$\\
\hline
$10000$ & $4760$ & $20.5$ & $71900$ & $167$ & $58.1$ & $20.5$ & $3 \times 10^{-13}$ & $8.0 \times 10^{-15}$ & $3.74 \times 10^{-15}$\\
\hline
$20000$ & $-$ & $49.1$ & $-$ & $333$ & $-$ & $42.9$ & $-$ & $1.0 \times 10^{-14}$ & $-$\\
\hline
$50000$ & $-$ & $116$ & $-$ & $838$ & $-$ & $108$ & $-$ & $1.5 \times 10^{-14}$ & $-$\\
\hline
$100000$ & $-$ & $216$ & $-$ & $1680$ & $-$ & $213$ & $-$ & $1.8 \times 10^{-14}$ & $-$\\
\hline
$200000$ & $-$ & $441$ & $-$ & $3380$ & $-$ & $425$ & $-$ & $2.6 \times 10^{-14}$ & $-$\\
\hline
$500000$ & $-$ & $1330$ & $-$ & $8500$ & $-$ & $1070$ & $-$ & $3.4 \times 10^{-14}$ & $-$\\
\hline
$1000000$ & $-$ & $2700$ & $-$ & $17600$ & $-$ & $2330$ & $-$ & $3.9 \times 10^{-14}$ & $-$\\
\hline
\end{tabular}
\end{center}
\label{table_N_scaling}
\end{table}
\begin{itemize}
\item
Assembly time - Time taken to assemble the dense matrix versus extended sparse matrix.
\item
Factorization time - Time taken to factorize the dense matrix versus extended sparse matrix.
\item
Solve time - Time taken to solve the dense linear system versus the extended sparse linear system (once the factorization has been obtained).
\item
Error in residual - Comparision of $\Vert Ax-b\Vert_{\infty}$ and $\Vert A_{ex}x_{ex} - b_{ex}\Vert_{\infty}$.
\item
Error in log-det - Relative error of the log of the absolute value of the determinant of the dense matrix and the extended sparse matrix.
\end{itemize}

\begin{figure}[!htbp]
\subfigure[Assembly time versus system size]{
\includegraphics{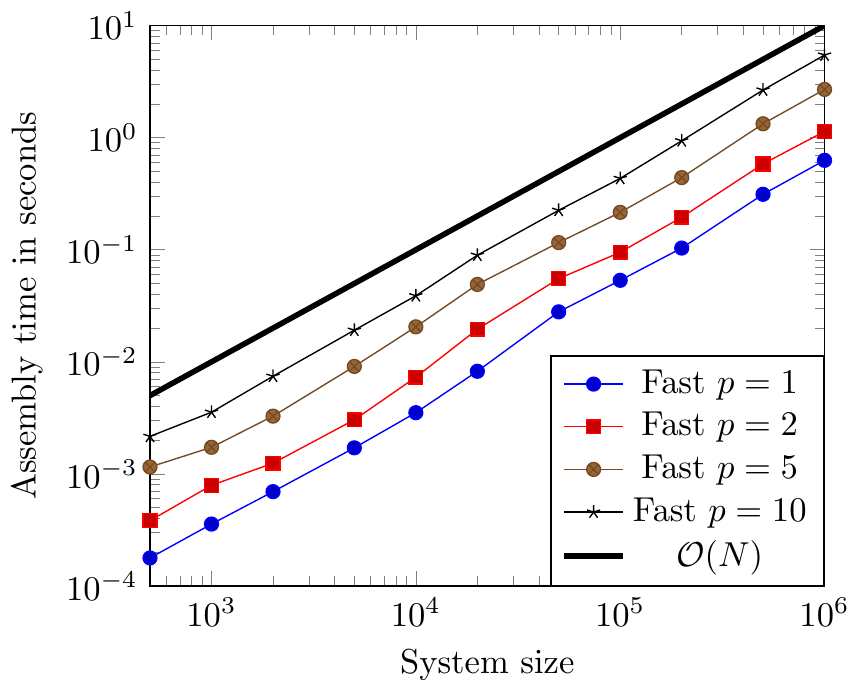}
}\subfigure[Factor time versus system size]{
\includegraphics{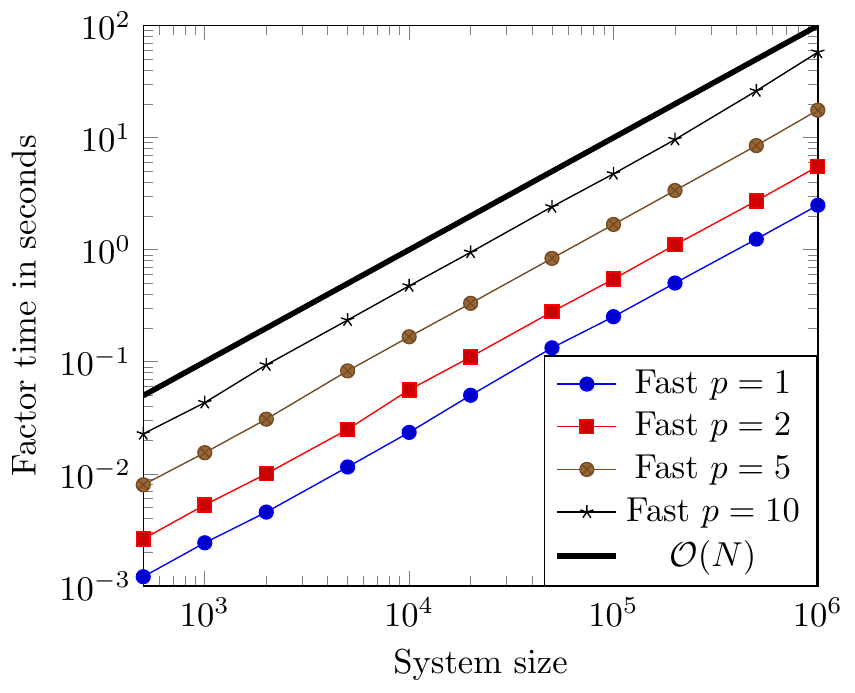}
}
\subfigure[Solve time versus system size]{
\includegraphics{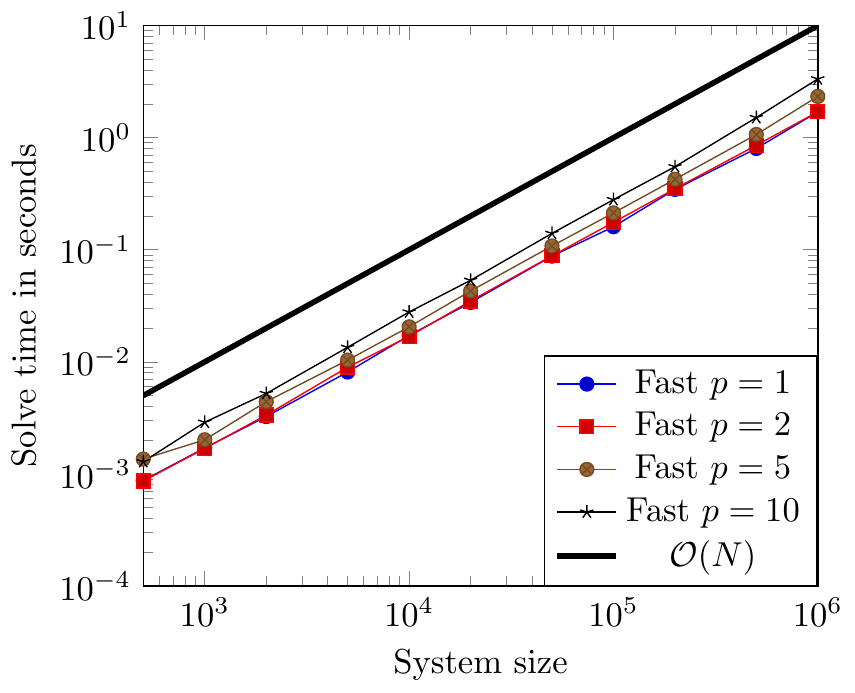}
}
\subfigure[Error versus system size]{
\includegraphics{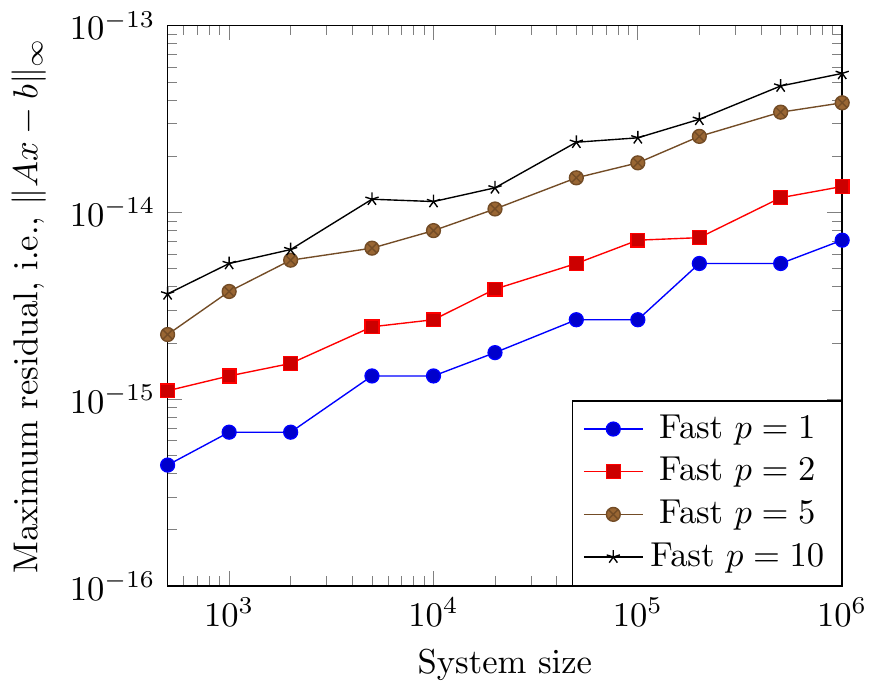}
}
\caption{Scaling of the algorithm with system size. From the benchmarks, it is clear that the computational cost for the fast algorithm scales as $\mathcal{O}(N)$ for assembly, factorization and solve stages, where $N$ is the number of unknowns. The maximum residual is less than $10^{-13}$ even for a system with million unknowns.}
\label{figure_benchmark1}
\end{figure}

Figure~\ref{figure_benchmark1} illustrates the scaling of the assembly, factorization and solve time with system size. The different components of the algorithm, i.e., assembly, factorization and solve, scale linearly in the number of unknowns. Also, as expected the pre-factor infront of the linear scaling increases with the semi-separable rank $p$, i.e., in our case the number of exponentials.

\subsection{Benchmark $2$}
In this benchmark, we illustrate the scaling of the time taken (assembly, factorization and solve) for algorithm with $p$, the number of exponentials added (equivalently the semi-separable rank).

\begin{figure}[!htbp]
\subfigure[Assembly time versus number of exponentials]{
\includegraphics{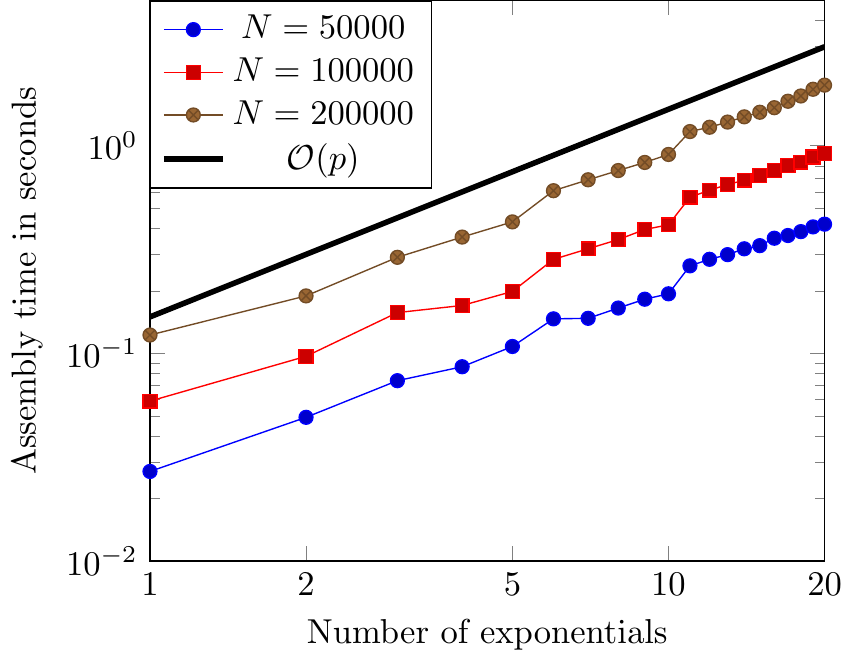}
}\subfigure[Factor time versus number of exponentials]{
\includegraphics{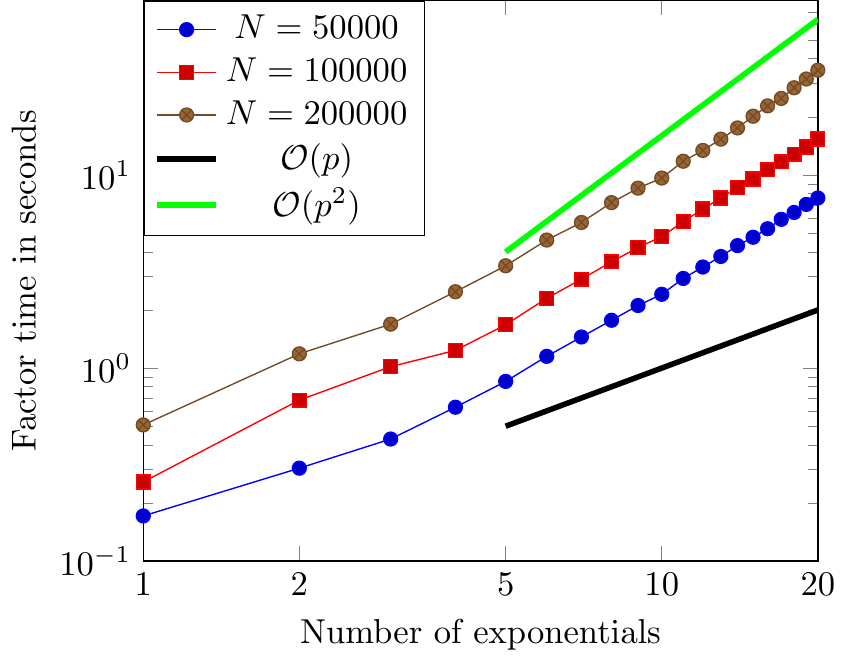}
}
\subfigure[Solve time versus number of exponentials]{
\includegraphics{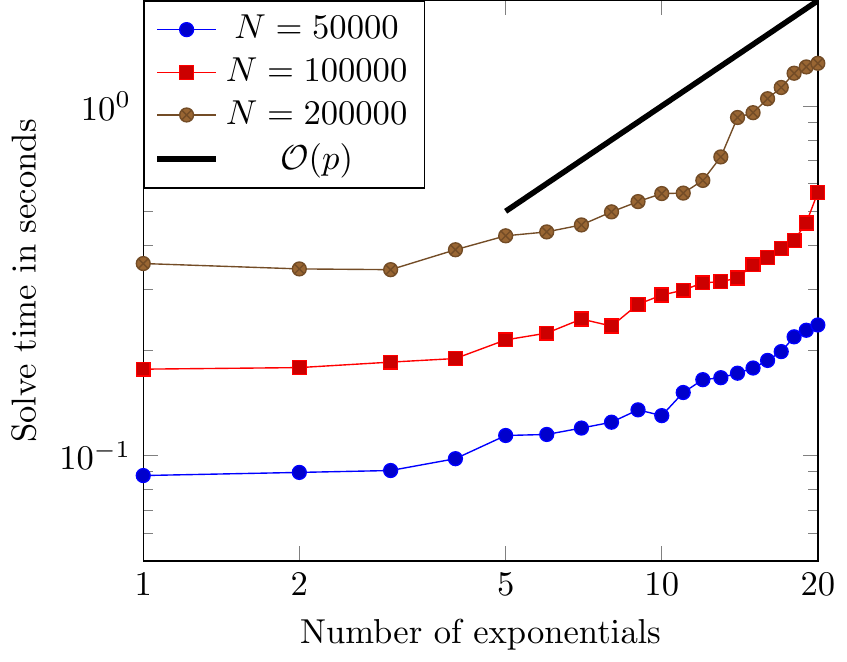}
}
\subfigure[Error versus number of exponentials]{
\includegraphics{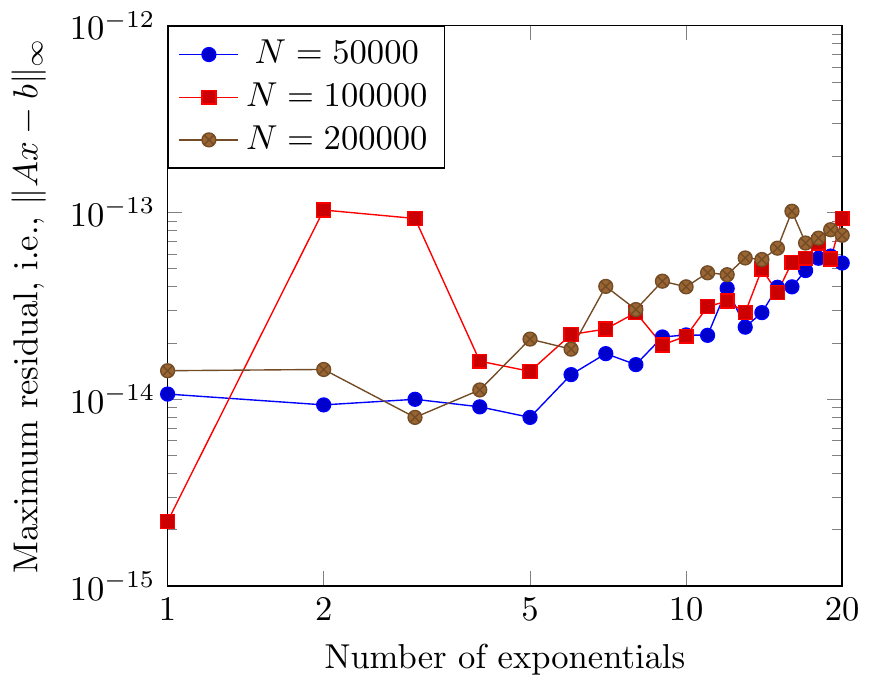}
}
\caption{Scaling of the fast algorithm with the number of exponentials added. From the benchmarks, it is clear that the computational cost for the fast algorithm scales as $\mathcal{O}(p)$ for assembly and solve stages, while it scales as $\mathcal{O}(p^2)$ for the factorization stage, where $p$ is the number of exponentials (equivalently the semi-separable rank). The maximum residual is less than $10^{-13}$ almost always.}
\label{figure_benchmark2}
\end{figure}

Figure~\ref{figure_benchmark2} illustrates the scaling of different parts of the algorithm  with the semi-separable rank $p$. Note that the assembly time scales linearly with the semi-separable rank, while the factorization time scales quadratically with the semi-separable rank as expected. The error in the solution seems to be more or less independent of the semi-separable rank.
\section{Conclusion}
The article discusses a numerically stable, generalized Rybicki Press algorithm, which relies on the fact that a semi-separable matrix can be embedded into a larger banded matrix. This enables $\mathcal{O}(N)$ inversion and determinant computation of covariance matrices, whose entries are sums of exponentials. This also immediately provides a fast matrix vector product for semi-separable matrices. This publication also serves to formally announce the release of the implementation of the extended sparse semi-separable factorization and the generalized Rybicki Press algorithm. The implementation is in C++ and is made available at \url{https://github.com/sivaramambikasaran/ESS}~\cite{ambikasaran2014ESS} under the license provided by New York University.

\FloatBarrier


\newpage

\end{document}